\documentclass[11pt, a4paper]{article}
\usepackage[margin = 3cm]{geometry}
\usepackage{amsthm, amsmath, amssymb}
\usepackage[numbers,sort&compress]{natbib}
\usepackage[colorlinks = true, citecolor = blue, linkcolor = blue, urlcolor = blue]{hyperref}
\usepackage{doi}
\usepackage{microtype}
\usepackage{booktabs}
\usepackage{tikz}
\usepackage{etoolbox}
\usepackage[affil-it]{authblk}

\usepackage[varbb]{newpx}
\everymath{\displaystyle}
\counterwithin{equation}{section}
\counterwithin{table}{section}
\counterwithin{figure}{section}
\newtheorem{thmcounter}{none}[section]
\theoremstyle{definition}

\theoremstyle{remark}
\newtheorem{remark}[thmcounter]{Remark}
\theoremstyle{plain}
\newtheorem{lemma}[thmcounter]{Lemma}
\newtheorem{proposition}[thmcounter]{Proposition}
\newtheorem{theorem}[thmcounter]{Theorem}

\newtheorem{conjecture}[thmcounter]{Conjecture}
\newtheorem{question}[thmcounter]{Question}

\newcommand{\Q}{\mathbb{Q}}

\newcommand{\C}{\mathbb{C}}
\newcommand{\fieldchar}{\operatorname{char}}

\newenvironment{acknowledgement}{\section*{Acknowledgements}}{}

\title{\bf Becker's Conjecture on Compositions of \\ Higher Degree Diagonal Forms}
\author[]{Hua-Lin Huang\thanks{Corresponding author: \texttt{hualin.huang@hqu.edu.cn}}}
\author[]{Yiming Liu}
\author[]{Huajun Lu}
\author[]{Yilun Tang}
\affil[]{Huaqiao University, Quanzhou, China}

\date{}

\begin{document}
    \maketitle
    \enlargethispage{4pt}

    \begin{abstract}
        \noindent Becker's conjecture asserts that a rational composition of sums of $d$-th powers of size $[r,s,n]$, with $d>2$, requires $n \ge rs$. We show that this formulation fails over every characteristic-zero field in odd degree, and in every degree over a suitable finite extension of $\Q$. 
        Over formally real fields in even degree, however, we prove that every polynomial composition satisfies $n\ge rs$.
    \end{abstract}

    \section{Introduction}

    \par For a field $F$, an integer $d\ge2$, and a tuple $X=(x_1,\ldots,x_r)$ of indeterminates, write
    \[
        S_d(X)=x_1^d+\cdots+x_r^d.
    \]
    A rational diagonal composition of degree $d$ and size $[r,s,n]$ is an identity
    \begin{equation}
        \label{equation:composition}
        S_d(X)S_d(Y)=\sum_{\nu=1}^n z_\nu(X,Y)^d,
        \qquad z_\nu\in F(X,Y) \setminus\{0\},
    \end{equation}
    where $X=(x_1,\ldots,x_r)$ and $Y=(y_1,\ldots,y_s)$ are disjoint tuples of algebraically independent variables. The products $x_iy_j$ always give a composition with $n=rs$.

    \par Shapiro \cite{Shapiro} records the following statement as Becker's conjecture.

    \begin{conjecture}
        \label{conjecture:becker}
        Suppose that $d>2$ and $\fieldchar F \nmid d!$. If the identity \eqref{equation:composition} holds, then $n\ge rs$.
    \end{conjecture}

    \par The conjecture was motivated by Becker's work \cite{Becker79, Becker82} on higher Pythagoras numbers. A related question for integer-coefficient polynomials was raised earlier by Nathanson \cite{Nathanson}. 

    \par The starting point for compositions of diagonal forms is the multiplication law for sums of squares. Hurwitz classified bilinear identities of size $[n,n,n]$, obtaining the celebrated ``$1,2,4,8$ Theorem" \cite{Hurwitz}. Pfister showed that rational compositions of size $[2^m,2^m,2^m]$ exist in degree two \cite{Pfister}. For compositions of higher degree diagonal forms, Gude \cite{Gude} obtained several nonexistence results. Very recently, some further results appear in \cite{Harrison}.

    \par We show that Becker's conjecture fails in every odd degree over $\Q$, and in every degree over a suitable finite extension of $\Q$. The construction uses a classical finite-difference identity to write a polynomial as a sum of $d$-th powers. Our second result proves the conjectured bound for polynomial summands over formally real fields in every even degree $d\ge4$.

    \par The note is organized as follows. As a preliminary, we recall a finite-difference identity in Section 2. In Section 3 we provide counterexamples to Becker's conjecture. In Section 4 we partially confirm the conjecture in the case of even degree over formally real fields. Finally, in Section 5, we suggest a modification of Becker's conjecture.

    \section{A finite-difference identity}

    \par All fields in the rest of the paper have characteristic zero. Let $f$ be a function of one variable. The difference operator $\Delta$ assigns to $f$ the function $\Delta f$ defined by
    \[
        \Delta f(T)=f(T+1)-f(T),
    \]
    and its iterates are defined by $\Delta ^{k+1} f(T) = \Delta ^{k} f(T+1) - \Delta ^{k} f(T)$. Induction on $k$ then yields the \emph{finite-difference formula}
    \begin{equation}\label{equation:difference}
        \Delta^{k}f(T)=\sum_{j=0}^{k}(-1)^{k-j}\binom{k}{j}f(T+j).
    \end{equation}

    \begin{lemma}
        \label{lemma:difference}
        For every integer $d\ge2$, the following identity holds:
        \begin{equation}\label{equation:centered}
            \sum_{i=0}^{d-1}(-1)^{d-1-i}\binom{d-1}{i}\left(T+i-\frac{d-1}{2}\right)^d=d!T.
        \end{equation}
    \end{lemma}
    
    \begin{proof}
        \par By the classical finite-difference formula for $d$-th powers,
        \[
            \Delta^{d-1}U^{d} = \sum_{j = 0}^{d - 1} (-1)^{d - 1 - j} \binom{d - 1}{j}(U + j)^d = \frac{d!}{2}(2U+d-1).
        \]
        For the second equality, by Proposition 1.9.2 in \cite{Stanley}, one can expand $U^d$ as
        \[
            U^{d}=\sum_{k=0}^{d}k!S(d,k)\binom{U}{k},
        \]
        where $S(d,k)$ is the Stirling number of the second kind, that is, the number of partitions of a $d$-element set into $k$ blocks, and $\tbinom{U}{k}=U(U-1)\cdots(U-k+1)/k!$ is the formal binomial coefficient, regarded as a polynomial in $U$. Since $\Delta\tbinom{U}{k}=\tbinom{U}{k-1}$, applying $\Delta^{d-1}$ term by term to the previous expression of $U^d$ gives
        \[
            \begin{aligned}
                \Delta^{d-1}U^{d}&=\sum_{k=d-1}^{d}k!S(d,k)\binom{U}{k-d+1}=(d-1)!S(d,d-1)+d!S(d,d)U \\
                &= d!(d-1)/2+d!U=d!(2U+d-1)/2.
            \end{aligned}
        \]
        It remains only to substitute $U=T-\frac{d-1}{2}$.
    \end{proof}

    \begin{proposition}
        \label{proposition:universal}
        \par Suppose that $(-1)^{d-1-i}\binom{d-1}{i}=\lambda_i^d$ with $\lambda_i\in F$, for every $0\le i\le d-1$. Then every polynomial $P$ over $F$ in any finite number of indeterminates satisfies
        \begin{equation}\label{equation:universalidentity}
            P=\sum_{i=0}^{d-1}\left[\lambda_i\left(\frac{P}{d!}+i-\frac{d-1}{2}\right)\right]^d.
        \end{equation}
    \end{proposition}

    \begin{proof}
        \par Substitute $T=P/d!$ into \eqref{equation:centered} and absorb each coefficient $(-1)^{d-1-i}\binom{d-1}{i}=\lambda_i^d$ into the corresponding $d$-th power.
    \end{proof}

    \section{Polynomial counterexamples}
    \par We apply the finite-difference identity to provide counterexamples to Becker's conjecture.
    \begin{theorem}\label{theorem:rootcase}
        \par Fix $d>2$ and a field $F$ containing $d$-th roots of the constants $(-1)^{d-1-i}\binom{d-1}{i}$, $0\le i\le d-1$. For every $r,s\ge1$, there is a polynomial composition of size $[r,s,d]$. Consequently the formulation in Conjecture \ref{conjecture:becker} fails whenever $rs>d$.
    \end{theorem}

    \begin{proof}
        \par Let $\lambda_i$ be a $d$-th root of the nonzero constant $(-1)^{d-1-i}\binom{d-1}{i}$, for $0\le i\le d-1$. Applying Proposition~\ref{proposition:universal} to the polynomial $P=S_d(X)S_d(Y)\in F[X,Y]$ gives
        \[
            S_d(X)S_d(Y)=z_1^d+\cdots+z_d^d,\qquad
            z_i=\lambda_{i-1}\left(\frac{S_d(X)S_d(Y)}{d!}+i-\frac{d+1}{2}\right),\quad 1\le i\le d.
        \]
        Each $z_i$ lies in $F[X,Y]$ and is nonzero, since $\lambda_{i-1}\neq0$ and the parenthesized factor has degree $2d>0$. Hence this is a polynomial composition of size $[r,s,d]$.
    \end{proof}
    For each $d$, the field
    \[
        F_d=\Q(\lambda_0,\ldots,\lambda_{d-1}),\qquad \lambda_i^d=(-1)^{d-1-i}\binom{d-1}{i},
    \]
    is a finite extension of $\Q$ on which the theorem applies. In particular, the theorem applies over $\overline\Q$ and $\C$. If $d$ is even, the constant at $i=0$ equals $(-1)^{d-1}=-1$, so $F_d$ is not formally real: $\lambda_0^d=-1$ makes $-1$ a square. The base-field extension is thus essential in even degree.

    \par For odd degrees there is a construction over $\Q$ itself, without adjoining roots.
    \begin{theorem}\label{theorem:odd}
        Let $d>2$ be odd. Over every characteristic-zero field and for every $r,s\ge1$, there is a polynomial composition of size
        \[
            [r,s,2^{d-1}].
        \]
        In particular the formulation in Conjecture \ref{conjecture:becker} over $\Q$ is false whenever $rs>2^{d-1}$.
    \end{theorem}

    \begin{proof}
        \par Put $P=S_d(X)S_d(Y)/d!\in F[X,Y]$. Since $d$ is odd, we have $(-1)^{d-1-i}=(-1)^{i}$ and $(-1)^{id}=(-1)^{i}$. Applying \eqref{equation:centered} to $P$ and absorbing each sign into the corresponding power, we obtain the decomposition
        \[
            \begin{aligned}
                S_d(X)S_d(Y) &= d!P
                =\sum_{i=0}^{d-1}(-1)^{i}\binom{d-1}{i}\left(P+i-\frac{d-1}{2}\right)^{d} \\
                &=\sum_{i = 0}^{d - 1} \sum_{j = 1}^{\binom{d - 1}{i}} z_{i,j}^{d},
            \end{aligned}
        \]
        where each $z_{i,j}=(-1)^{i}\left(P+i-\frac{d-1}{2}\right)$. This is a polynomial composition of size $[r,s,2^{d-1}]$.
    \end{proof}

    \begin{remark}
        \par These counts are upper bounds for the number of summands, not assertions of minimality. Repeated summands are permitted in \eqref{equation:composition}. All coefficients lie in $\Q$, but need not be integers. Thus these examples do not settle integer-coefficient versions of the problem, such as the question studied by Nathanson \cite{Nathanson}.
    \end{remark}

    \section{Even degrees over formally real fields}

    \par Note that the alternating coefficients in \eqref{equation:centered} cannot be absorbed into even powers over an ordered field. In fact, we have a stronger obstruction in the polynomial case.

    \begin{theorem}\label{theorem:positive}
        \par Let $F$ be a formally real field and let $d\ge4$ be even. Suppose that
        \begin{equation}\label{equation:polynomial}
            S_d(X)S_d(Y)=\sum_{\nu=1}^n z_\nu^d,
            \qquad z_\nu\in F[X,Y]\setminus \{0\}.
        \end{equation}
        Then each $z_\nu$ is of the form
        \[
            z_\nu=a_\nu x_{i(\nu)}y_{j(\nu)},\qquad a_\nu\in F^\times.
        \]
        In particular $n\ge rs$, and this bound is attained.
    \end{theorem}
    
    \begin{proof}
        \par Recall that a formally real field is orderable and has characteristic zero. Consequently a sum of squares of polynomials can vanish only if each summand vanishes. The same holds for sums of $d$-th powers, since $d$ is even.

        \par We first show that the $z_\nu$ are bilinear. Let $D=\max_\nu\deg_X z_\nu$, where $\deg_X$ denotes the total degree in $X$ alone. If $D>1$, let $H_\nu$ be the $X$-homogeneous part of degree $D$ of $z_\nu$. Comparing $X$-degrees $dD$ in \eqref{equation:polynomial} gives $\sum_\nu H_\nu^d=0$, so all $H_\nu$ vanish, contradicting the choice of $D$. Thus $\deg_X z_\nu\le1$, and similarly $\deg_Y z_\nu\le1$. Setting $X=0$ and $Y=0$ in turn shows that $z_\nu(0,Y)=z_\nu(X,0)=0$. These facts exclude constants and terms involving only one tuple, giving
        \[
            z_\nu=\sum_{i=1}^r x_i L_{\nu i}(Y),
        \]
        with homogeneous linear forms $L_{\nu i}$.

        \par Next, we claim that for each $\nu$ at most one of the forms $L_{\nu i}(Y)$ is nonzero. Suppose that $L_{\nu i}(Y)\neq0$ and $L_{\nu i'}(Y)\neq0$ for some $i\ne i'$. Then the coefficient of $x_i^{d-2}x_{i'}^2$ on the left side of \eqref{equation:polynomial} is zero, since $d\ge4$. On the right it is $\binom d2\sum_\nu L_{\nu i}^{d-2}L_{\nu i'}^2$. Dividing by $\binom d2$ therefore yields
        \begin{equation}
            \label{equation:mixed-coefficient}
            \sum_\nu\left(L_{\nu i}^{(d-2)/2}L_{\nu i'}\right)^2=0.
        \end{equation}
        Hence $L_{\nu i}^{(d-2)/2}L_{\nu i'}=0$ for every $\nu$. Since the polynomial ring has no zero divisors, we get a contradiction. Applying the same argument to the $y_j$, we conclude that every nonzero summand has the asserted form $z_\nu=a_\nu x_{i(\nu)}y_{j(\nu)}$.

        \par Finally, for every pair $(i,j)$ with $1\le i\le r$ and $1\le j\le s$, the monomial $x_i^dy_j^d$ has coefficient $1$ on the left of \eqref{equation:polynomial}, so at least one index $\nu$ satisfies $(i(\nu),j(\nu))=(i,j)$. This yields that $n \ge rs$.

        \par Taking $z_{ij}=x_iy_j$ for all pairs achieves equality.
    \end{proof}

    \begin{remark}
        \par The exclusion of $d=2$ is essential: the coefficient $x_i^{d-2}x_j^2$ would then be a pure power already present on the left. 
    \end{remark}

    \section{The remaining rational problem}

    \par The preceding results suggest restricting the arithmetic setting rather than imposing a universal inequality over all characteristic-zero fields.
    \begin{question}\label{question:rational}
        Let $F$ be formally real and let $d\ge4$ be even. Does every rational composition \eqref{equation:composition} satisfy $n\ge rs$?
    \end{question}
    Theorem~\ref{theorem:positive} settles the polynomial case, but its proof does not extend by simply clearing denominators: a common denominator introduces an additional factor on the left and invalidates the degree argument. 

    \par Gude's results, as reported in \cite{Shapiro}, give affirmative cases over $\Q$, including $d=4$ and $rs\le15$. Another restriction is to require the summands to be bilinear, or linear in $Y$ over $F(X)$ as in \cite{Harrison}. These conditions exclude the finite-difference constructions, whose summands are nonhomogeneous and have much larger degree. They lead to separate composition problems and should not be conflated with unrestricted rational identities.

    \par The bounds $d$ and $2^{d-1}$ obtained above are not claimed to be minimal. Determining the smallest number of summands for given $F,r,s,d$, and determining which finite constant extensions admit the $d$-term construction, remain open questions.

    \begin{acknowledgement}
        \par H.-L. Huang sincerely thanks Professor Zach Teitler for his insightful questions on \cite{Harrison} that motivated this short note. This research is partially supported by the Key Program of the Natural Science Foundation of Fujian Province (Grant No. 2024J02018) and the National Natural Science Foundation of China (Grant No. 12371037).
    \end{acknowledgement}

    \bibliography{references}

@article{Becker79,
  author  = {Becker, Eberhard},
  title   = {{Summen $n$-ter Potenzen in K\"orpern}},
  journal = {Journal f\"ur die reine und angewandte Mathematik},
  volume  = {307/308},
  pages   = {8--30},
  year    = {1979},
  doi     = {10.1515/crll.1979.307-308.8}
}

@incollection{Becker82,
  author    = {Becker, Eberhard},
  title     = {{The real holomorphy ring and sums of $2n$-th powers}},
  booktitle = {{G\'eom\'etrie alg\'ebrique r\'eelle et formes quadratiques}},
  editor    = {Colliot-Th\'el\`ene, Jean-Louis and Coste, Michel and Mah\'e, Louis},
  series    = {Lecture Notes in Mathematics},
  volume    = {959},
  pages     = {139--181},
  publisher = {Springer},
  address   = {Berlin},
  year      = {1982},
  doi       = {10.1007/BFb0062253}
}

@phdthesis{Gude,
  author  = {Gude, Ulf-Gerhard},
  title   = {{\"Uber die Nichtexistenz rationaler Kompositionsformeln bei Formen h\"oheren Grades}},
  school  = {Universit\"at Dortmund},
  address = {Dortmund},
  year    = {1988}
}

@article{Harrison,
  author  = {Huang, Hua-Lin and Liao, Lili and Lu, Huajun and Ye, Yu and Zhang, Chi},
  title   = {{Harrison Center and Products of Sums of Powers}},
  journal = {Communications in Mathematics and Statistics},
  volume  = {14},
  number  = {1},
  pages   = {1--8},
  year    = {2026},
  doi     = {10.1007/s40304-023-00367-1}
}

@article{Hurwitz,
  author  = {Hurwitz, Adolf},
  title   = {{\"Uber die Komposition der quadratischen Formen von beliebig vielen Variabeln}},
  journal = {Nachrichten von der Gesellschaft der Wissenschaften zu G\"ottingen, Mathematisch-Physikalische Klasse},
  pages   = {309--316},
  year    = {1898}
}

@article{Nathanson,
  author  = {Nathanson, Melvyn B.},
  title   = {{Products of Sums of Powers}},
  journal = {Mathematics Magazine},
  volume  = {48},
  pages   = {112--113},
  year    = {1975},
  doi     = {10.1080/0025570X.1975.11976458}
}

@article{Pfister,
  author  = {Pfister, Albrecht},
  title   = {{Multiplikative quadratische Formen}},
  journal = {Archiv der Mathematik},
  volume  = {16},
  pages   = {363--370},
  year    = {1965},
  doi     = {10.1007/BF01220043}
}

@book{Shapiro,
  author    = {Shapiro, Daniel B.},
  title     = {{Compositions of Quadratic Forms}},
  series    = {de Gruyter Expositions in Mathematics},
  volume    = {33},
  publisher = {Walter de Gruyter},
  address   = {Berlin},
  year      = {2000},
  doi       = {10.1515/9783110824834}
}

@book{Stanley,
  author    = {Stanley, Richard P.},
  title     = {{Enumerative Combinatorics, Volume~1}},
  edition   = {2nd},
  series    = {Cambridge Studies in Advanced Mathematics},
  number    = {49},
  publisher = {Cambridge University Press},
  address   = {Cambridge},
  year      = {2012},
  doi       = {10.1017/CBO9781139058520}
}
\end{document}